# Causal inference in longitudinal studies with history-restricted marginal structural models

## Romain Neugebauer and Mark J. van der Laan

*Division of Biostatistics, School of Public Health, University of California, Berkeley*
*e-mail:* `romain@stat.berkeley.edu`; `laan@stat.berkeley.edu`

## Marshall M. Joffe

*Department of Biostatistics and Epidemiology, University of Pennsylvania,*
*School of Medicine e-mail:* `mjoffe@mail.med.upenn.edu`

## Ira B. Tager

*Division of Epidemiology, School of Public Health, University of California, Berkeley*
*e-mail:* `ibt@berkeley.edu`

**Abstract:** A new class of Marginal Structural Models (MSMs), History-Restricted MSMs (HRMSMs), was recently introduced for longitudinal data for the purpose of defining causal parameters which may often be better suited for public health research or at least more practicable than MSMs (6, 2). HRMSMs allow investigators to analyze the causal effect of a treatment on an outcome based on a fixed, shorter and user-specified history of exposure compared to MSMs. By default, the latter represent the treatment causal effect of interest based on a treatment history defined by the treatments assigned between the study's start and outcome collection. We lay out in this article the formal statistical framework behind HRMSMs. Beyond allowing a more flexible causal analysis, HRMSMs improve computational tractability and mitigate statistical power concerns when designing longitudinal studies. We also develop three consistent estimators of HRMSM parameters under sufficient model assumptions: the Inverse Probability of Treatment Weighted (IPTW), G-computation and Double Robust (DR) estimators. In addition, we show that the assumptions commonly adopted for identification and consistent estimation of MSM parameters (existence of counterfactuals, consistency, time-ordering and sequential randomization assumptions) also lead to identification and consistent estimation of HRMSM parameters.

## 1. Motivation

Longitudinal epidemiological studies are increasingly becoming more interested in the time-dependent effects of various exposures on human health outcomes.

This is particularly true for studies concerned with the effects of chronic exposure to ambient air pollutants. Cohort studies with multiple time-specific estimates of exposure have been emphasized by the U.S. Environmental Protection Agency (EPA) as the preferred study design to address these issues (1). Examples of such concerns can be found for other health relevant exposures, such as time-specific patterns of physical activity on cardiovascular outcomes and obesity (22).

Typically, cohort studies collect data at regular time intervals for all members of the cohort. In practice, each collection time represents a "window" of time over which data are collected. Information on the exposure of interest, also referred to as "treatment of interest", and other relevant covariates are obtained for the interval between successive data collection time points.

Currently, the principal tools used by epidemiologists for the analysis of cohort data are conditional, association models (e.g., logistic, pooled logistic, Cox proportional hazards). The time-dependence of exposure effects usually are addressed in one of two ways: 1) ignored, in that only baseline exposure is considered or the time-dependent exposures are summarized in one cumulative measure whose effect is assumed to be confounded by baseline variables only (e.g. "pack-years" is a cumulative measure of exposure to cigarette smoke commonly used in occupational epidemiology) ; and 2) risk-set sampling, in which the exposure measure (possibly cumulative as just described) is updated at specific time points and exposure effects estimated based on the probability of outcome within groups classified by exposure in the current or most recent time interval. Thus, the full complexity of the effect of an exposure history on a specific health outcome is typically lost in conditional, association models. More importantly, these models also suffer from a major drawback: they typically produce biased effect estimates in scenarios involving time-dependent confounders also "on one of the causal pathways" of interest (16, 19). Finally, these models do not provide direct population-level estimates of exposure effects, which often are of most relevance to public health. Only conditional associations between the exposure and outcome can be directly derived from traditional regression models. Indirectly however, causal inferences can be drawn from such models in particular cases and under additional assumptions. In such scenarios, it is important to note that only conditional (i.e. not population-level) effect estimates can typically be derived from such models.

MSMs define parameters with a direct population-level causal interpretation. Moreover, estimators of MSM parameters can account properly for time-dependent confounders also on one of the causal pathways of interest between the exposure and health outcome (19). Thus, these models can represent accurately population-level effects of histories of exposures on the health outcomes: the cumulating effects of exposures as well as acute effects from exposures experienced over shorter time intervals. This is particularly important for exposures, such as ambient air pollutants, whose acute and chronic exposure effects lead to different health outcomes and/or contribute together in the occurrence of serious health outcomes such as heart attacks and death due to diseases of the heart and lung (1).

The current MSM methodology represents exposure effects by considering the exposure histories over *the entire* follow-up intervals that precede the occurrence of time-specific health outcomes. In cases where consideration of the entire time interval prior to an event occurrence may not be relevant, based on subject matter, this omnibus handling of time possess an important limitation. This paper presents an extension of current MSM methodology to allow for more flexible analysis of time effects of exposure, based on *a priori* or *ad hoc* considerations of specific periods of time antecedent to an event.

Note that this modeling extension was originally proposed and motivated in earlier work with discretized failure-time data and called partially marginal structural models (6, 2). We develop here the formal statistical framework behind these models, compare it to the MSM framework and further motivate this approach with two other research problems.

First, we undertook a study to determine the extent to which reductions in ambient air pollution consequent to regulations propagated since 1980 by the California Air Resources Board to reduce air pollution in the Los Angeles (LA) Basin lead to measurable health benefits. The basic time unit for the data was the quarter (3 months), and we had 84 such time units. The geographic area of interest was divided into 195 10 x 10 km grids, based on know patterns of air pollutants and meteorology in the LA Basin. The exposures of interest were quarterly concentrations of ambient air pollutants (e.g., ozone, oxides of nitrogen, particulate matter with a mass medium aerodynamic diameter of 10 microns or less). A variety of health outcomes are being considered: quarterly hospital discharges and mortality rates for various chronic lung and cardiovascular diseases. Population denominators and demographic data are available on a quarter-spatial unit-specific basis. Over the 20 years encompassed by the study, there have been large temporal trends in demography that have lead to changes in population susceptibility to certain diseases of interest. For example, there has been a large influx of Mexicans into the study area. Mexicans are known to have decreased risks for asthma and increased risks for diabetes mellitus, an important underlying risk factor for cardiovascular disease. Moreover, since many of these immigrants are of low socioeconomic status, they may be more likely to live in closer proximity to sources of ambient pollutants (i.e., near Freeways). Thus, demography is an important temporal confounder. Changes in medical care over the study period also has affected the occurrence of health outcomes. In our descriptive analyses, important temporal trends for hospital discharges for asthma, chronic obstructive lung disease and various cardiovascular disease were observed. Temporal trends of disease-specific mortality are expected, e.g. there has been a decline in age-specific death rates from specific heart disease due to improvements in medical care. Based on the above, "time" becomes an important confounding variable to capture all of the unmeasured temporally-related factors that we have not measured and the residual temporal confounding for those factors that we have measured. Since our initial analysis focused on hospital discharges for asthma in children ages birth to 19 years, a central issue that emerged was the relevant exposure time for investigation of the exposure effect on the outcome rates. Based on available data on the ef-fects of ambient air pollution on hospital discharges for asthma, it did not seem reasonable to extend the exposure period much beyond the 12 months prior to a given quarter. Consequently, there was a need to modify the how exposure histories are handled in the existing MSM methodology.

Second, the proposed HRMSM methodology also has application to panel data that are being collected as part of a study of the relation between responses to short-term increases in ambient air pollutants and the long-term changes in symptoms and disease severity in children with asthma (Fresno Asthmatic Children's Environment Study - FACES). In this study, subjects participate in multiple 14-day panels during which time each subjects provides daily data on lung function, respiratory symptoms and daily activities. Analyses focus on the causal relation between daily symptoms or lung function and exposure to one or another pollutant over one or more days prior outcome report. Confounders for these analyses relate to meteorology and the effects of other pollutants not of primary interest in a given analysis. Virtually all studies to date evaluate the temporal relation between pollutants and symptoms/lung function through one of several modeling techniques based on traditional regression models where the study of the exposure effect is limited to: 1) specific lag effects (e.g., exposure on one or more days prior to outcome report) ; 2) specific average-over-several-days effects, or 3) polynomial-distributed-lag-functions effects (5). These approaches underscore the implicit understanding that the exposure effects of interest are limited in time.

To reflect these subject-matter considerations, we propose HRMSMs as an extension of MSMs to provide greater flexibility in the evaluation of exposure effects over time within a rigorous statistical framework for causal inference based on counterfactuals. Once the investigator specifies the time frame over which the effect of the pollutant exposure is of interest, a recently developed data-adaptive methodology for model selection (26) can provide guidance on how each level of the pollutant during the specified time frame should enter the HRMSM (e.g. moving averages or lags) that represents the hypothesized causal relationship between the pollutant and asthma outcome.

In section 2, we introduce HRMSMs and in section 3 present three estimators of HRMSM parameters developed based on a so called $t$-specific counterfactual framework. HRMSMs are then compared to MSMs in section 4 before illustrating the application of the proposed HRMSM methodology with a real-life application in section 5. Finally, we summarize and discuss the results presented in this article in section 6. For completeness, we include technical details in a comprehensive appendix.

## 2. History-Restricted Marginal Structural Model

### 2.1. Data structure and question of interest

For all experimental units in a random population sample of size $n$, we observe a treatment regimen $(A(0), \ldots, A(K))$ over time $t = 0, \ldots, K$ and a covariate process $(L(0), \ldots, L(K+1))$ measured at baseline and after a new treatment is assigned. The covariate $L(t)$ is measured after $A(t-1)$ and before $A(t)$. Note that $K+1$ represents the length of the treatment regimen in the appropriate unit of time and $n$ the sample size.

In the formal counterfactual framework for longitudinal studies (26), the data are represented as $n$ independent and identically distributed (i.i.d) realizations of:

$$O = (L(0), A(0), L(1), A(1), \ldots, L(K), A(K), L(K+1)) = (\bar{A}(K), \bar{L}(K+1)) \sim P,$$

where $P$ represents the distribution of the stochastic process $O$, referred to as the observed data, and the general notation $\bar{\cdot}(t)$ represents the history of the variable '·' between time 0 and time $t$: a) $\bar{\cdot}(t) = (\cdot(0), \ldots, \cdot(t))$ if $t \geq 0$ and b) $\bar{\cdot}(t)$ is empty if $t < 0$. We extend this notation with the notation $\bar{\cdot}(t_-, t_+)$ to represent the history of the variable $\cdot$ between time points $t_-$ and $t_+$: where a) $\bar{\cdot}(t_-, t_+) = (\cdot(t_-), \ldots, \cdot(t_+))$ if $t_- \leq t_+$, and b) $\bar{\cdot}(t_-, t_+)$ is empty otherwise. We thus have $\bar{\cdot}(t) = \bar{\cdot}(0, t)$.

We denote with $s$ an integer between 0 and $K+1$ specified by the investigators to represent the history size for which the exposure effect is of interest in the analysis. For $s - 1 \leq t \leq K$, we define $V(t - s + 1)$ as a subset of $(\bar{A}(t-s), \bar{L}(t-s+1))$ to represent time-specific baseline variables, i.e. $V(t-s+1) \subset (\bar{A}(t-s), \bar{L}(t-s+1))$ and $V(t-s+1)$ represents the baseline variables if $t - s + 1$ is considered the study's starting time, i.e. the time from when the treatments and outcomes are considered of interest. In particular, we define a subset of the 'true' baseline variables as $V \equiv V(0) \subset L(0)$. We denote the time-dependent outcome with $Y(t)$: $Y(t) \in L(t)$, and $\mathcal{T}$ represents the set of time points $t$ such that the outcome, $Y(t+1)$, is of interest. We have $\mathcal{T} \subset \{0, \ldots, K\}$. Typically $\mathcal{T} = \{0, \ldots, K\}$ except when one is interested in the outcome collected at the end of the study only, i.e. when $\mathcal{T} = \{K\}$.

The question of interest is to investigate the causal effect of the treatment on the time-dependent outcomes of interest. In the literature, this problem has been addressed with MSMs. We propose in this article to address the same problem with the proposed HRMSMs. In MSM-based causal inference, the investigation of the causal relationship of interest relies on a representation of the effects of the treatment history between baseline and time point $t$, $\bar{A}(t)$, on the time-dependent outcome, $Y(t+1)$, for all $t \in \mathcal{T}$ (see figure 1). In HRMSM-based causal inference however, the investigation of the causal relationship of interest relies on a representation of different causal effects: the effects of the treatment history between time points $t - s + 1$ and $t$, $\bar{A}(t-s+1, t)$, on the time-dependent outcome, $Y(t+1)$, for all $t \in \mathcal{T}$. Compared to MSM-based causal inference, the effect of the treatment is thus investigated for a history of treatment that is restricted by the investigators based on considerations discussed later in this article (see figure 2). In other words, MSMs and HRMSMs can be viewed as two different statistical alternatives for the investigation of any given causal relationship, each providing different information and description of the causal effect of interest. We argue that an HRMSM-based causal inference strategy may often be more suitable or at least more practicable than an MSM-based causal inference strategy for public health research.

### 2.2. Assumptions

**Existence of counterfactuals:** We assume the existence of the following treatment-specific processes, also referred to as a counterfactual processes (11, 20), $\bar{L}_{\bar{a}(K)}(K + 1)$ for every treatment regimen $\bar{a}(K) = (a(0), \ldots, a(K)) \in \mathcal{A}_V(K)$, where $\mathcal{A}_V(K)$ designates all possible treatment regimens between time points 0 and $K$ as a function of the baseline variable $V$ only. In other words, $\mathcal{A}_V(K)$ is the support of the conditional distribution of $\bar{A}(K)$ given $V$, $g(\bar{A}(K) \mid V)$. We denote the so-called full data process with $X$ and and its distribution with $F_X$: $X = \left( V, \left( \bar{L}_{\bar{a}(K)}(K + 1) \right)_{\bar{a}(K) \in \mathcal{A}_V(K)} \right) \sim F_X$.

Note that the existence of the counterfactual process $\bar{L}_{\bar{a}(K)}(K + 1)$ for every treatment regimen $\bar{a}(K) \in \mathcal{A}_V(K)$ implies the existence of the counterfactual processes $\bar{L}_{\bar{a}(t)}(t + 1) \equiv \bar{L}_{\bar{a}(t), \bar{A}(t+1, K)}(t + 1) \subset X$ for every $t = 0, \ldots, K - 1$ and every treatment regimen $\bar{a}(t) = (a(0), \ldots, a(t)) \in \mathcal{A}_V(t)$ where $\mathcal{A}_V(t)$ designates all possible treatment regimens between time points 0 and $t$, i.e. the support of the conditional distribution of $\bar{A}(t)$ given $V$, $g(\bar{A}(t) \mid V)$. We have $\mathcal{A}_V(t) = \{\bar{a}(t) : \exists \ \bar{a}'(K) \in \mathcal{A}_V(K) \ \ \bar{a}(t) = \bar{a}'(t)\}$ for $t = 0, \ldots, K - 1$ and $\mathcal{A}_V(t)$ is thus entirely defined by $\mathcal{A}_V(K)$.

**Consistency assumption:** At any time point $t$, we assume the following link between the observed data and the counterfactuals: $L(t) = L_{\bar{A}(K)}(t)$. Under this assumption, we have: $O = (\bar{A}(K), \bar{L}_{\bar{A}(K)}(K + 1)) \equiv \phi(\bar{A}(K), X)$, where $\phi$ is a specified function of the full data process $X$. This notation indicates that the problem can be treated as a missing data problem. Only the treatment-specific process associated with the observed treatment $\bar{A}(K)$ is observed; the others are missing.

**Temporal Ordering assumption:** At any time point $t$, we assume that any treatment-specific variable can only be affected by past treatments: $L_{\bar{a}(K)}(t) = L_{\bar{a}(t-1)}(t)$ for $t = 0, \ldots, K + 1$, where $L_{\bar{a}(-1)}(0) = L(0)$. This assumption is typically implied by the data collection procedure: the covariate $L(t)$ (which contains the outcome at time $t$) is measured after $A(t - 1)$ and before $A(t)$.

**Sequential Randomization Assumption (SRA):** At any time point $t$, we assume that the observed treatment is independent of the full data, given the data observed up to time point $t$: $A(t) \perp X \mid \bar{A}(t-1), \bar{L}(t)$. Under the SRA, the treatment mechanism, i.e. the conditional density or probability of $\bar{A}(K)$ given $X$: $g(\bar{A}(K) \mid X)$, becomes:

$$g(\bar{A}(K) \mid X) = \prod_{t=0}^{K} g(A(t) \mid \bar{A}(t-1), X) \overset{SRA}{=} \prod_{t=0}^{K} g(A(t) \mid \bar{A}(t-1), \bar{L}(t)).$$

The SRA implies coarsening at random (4) and, thus, the likelihood of the observed data factorizes into two parts: a so-called $F_X$ and $g$ part. The $F_X$ part of the likelihood only depends on the full data process distribution, and the $g$ part of the likelihood only depends on the treatment mechanism. As a consequence of this factorization of the likelihood under the SRA, we now denote the distribution of the observed data with $P_{F_X,g}$ and the likelihood of $O$ is:

$$\mathcal{L}(O) \stackrel{SRA}{=} \overbrace{f(L(0)) \underbrace{\prod_{t=1}^{K+1} f(L(t) \mid \bar{L}(t-1), \bar{A}(t-1))}_{Q_{F_X}}}^{F_X \text{ part}} \overbrace{g(\bar{A}(K) \mid X)}^{g \text{ part}}.$$

In addition, we denote the set of conditional densities or probabilities that define the $F_X$ part of the likelihood, except for $f(L(0))$ with $Q_{F_X}$.

### *2.3. HRMSM and causal parameter of interest*

We define an HRMSM as a model for a feature (e.g. expectation) of the distribution of the counterfactual outcomes, $Y_{\bar{A}(t-s),\bar{a}(t-s+1,t)}(t+1)$, corresponding with treatment interventions between time point $t-s+1$ and $t$ only (treatments between time points 0 and $t-s$, $\bar{A}(t-s)$, are left random), possibly conditional on the baseline covariates at time $t-s+1$, $V(t-s+1)$, for all possible treatment interventions $\bar{a}(t-s+1,t)$ and $t \in \mathcal{T}_s$ where $\mathcal{T}_s$ represents the set of time points $t$ such that the outcome $Y(t+1)$ is of interest and $t \geq s-1$, $\mathcal{T}_s = \{t \in \mathcal{T} : t \geq s-1\}$. We have $\mathcal{T}_s \subset \{s-1, \ldots, K\}$. Typically we will have $\mathcal{T}_s = \{s-1, \ldots, K\}$.

In section 4, we discuss the interpretation of HRMSM parameters and how they represent the causal relationship of interest for a given value for $s$ compared to MSM parameters. By convention in this article, the random portion of the treatment history defining counterfactuals is excluded from the counterfactual notation and thus we adopt the following notations: $Y_{\bar{A}(t-s),\bar{a}(t-s+1,t)}(t+1) \equiv Y_{\bar{a}(t-s+1,t)}(t+1)$.

Typically and specifically in this article, one is interested in average causal effects per stratum $V(t-s+1)$ of the population which can be represented by causal parameters defined by HRMSMs for $E_{F_X,g}(Y_{\bar{a}(t-s+1,t)}(t+1) \mid V(t-s+1))$ for $t \in \mathcal{T}_s$. Similar to MSM analyses, investigators can choose to model $E_{F_X,g}(Y_{\bar{a}(t-s+1,t)}(t+1) \mid V(t-s+1))$ for each $t \in \mathcal{T}_s$ separately with a stratified model or simultaneously with a pooled model (8). We denote the causal parameters defined by a stratified and pooled HRMSM with $\beta_t(F_X, g \mid \cdot)$ (one parameter for each $t$) and $\beta(F_X, g \mid \cdot)$ (single parameter for all $t$), respectively, to indicate that they are mappings from the space of distributions $(F_X, g)$ to the space of real numbers and that these mappings are functions of modelling assumptions represented by '·'. Similar to the two MSM approaches developed in the literature (10, 8), a parametric or nonparametric HRMSM approach can be adopted, each relying on different modelling assumptions.

Note that unlike the class of MSMs, HRMSMs are introduced as a class of mixed full and observed data models, since an HRMSM models the marginal distribution of counterfactuals, where part of the treatment history is left random. The distribution of the random portion of the treatment history is defined by the treatment mechanism, $g$, and, thus, is identified by the observed data.

## 3. HRMSM estimation

For clarity, we only present in this section a summary of the theoretical results detailed in the appendix as they relate to HRMSM estimation. For conciseness, we illustrate their practical applications with the example of one estimator only.

### 3.1. The t-specific counterfactual framework

The time-specific ($t$-specific) counterfactual framework can be viewed as an extension of the conventional counterfactual framework on which MSM-based causal inference relies. This latter mathematical construct was described in sections 2.1 and 2.2. It provides the rigorous framework to define, identify and estimate MSM parameters with the full and observed data based on a sufficient set of assumptions presented in section 2.2. We describe the $t$-specific counterfactual framework in detail in appendix A and introduce it because it allows us to generalize the MSM estimation procedures to HRMSM estimation procedures with minimum effort. After linking the $t$-specific counterfactual framework to the conventional counterfactual framework in appendix A, we show in appendix B that the sufficient set of assumptions for MSM estimation presented in section 2.2 is also sufficient for HRMSM estimation. An important consequence is that the choice of HRMSM-based causal analysis over MSM-based causal analysis becomes only a matter of practical considerations (statistical power and computing issues) and, above all, subject-matter considerations (the relevance of the causal effect representation to public health research). These considerations are developed in section 4.

### 3.2. Link between the conventional and t-specific counterfactual frameworks

Figure 5 in the appendix illustrates based on an example of a longitudinal study with short follow-up the link between the longitudinal data representation in the conventional counterfactual framework and its representation in the time-specific counterfactual framework. Note that in the conventional counterfactual framework the data are approached as a single entity, $O$, in the sense that the treatment is defined once and for all as a history $\bar{A}(K)$ and the outcome is time-dependent, $Y(t) \in L(t)$. On the other hand, in the $t$-specific counterfactual framework the data are viewed as layers of separate entities, $O^t$. For each $O^t$, the treatment and outcome of interest are redefined along with the baseline covariates (highlighted in yellow/light gray on figure 5). Note that for each $O^t$, the outcome is no longer time-dependent but correspond with the last outcome collected, $Y^t \in L^t(t+1)$, (highlighted in orange/dark gray on figure 5) and the treatment history size is fixed to a user-specified value $s = 2$. In the conventional counterfactual framework, the investigator examines the effect of $\bar{A}(K)$ on $Y(t)$ for all $t \in \mathcal{T}$ based on MSMs for the full data associated with $O$ whereas in the $t$-specific counterfactual framework, the investigator can examine the effect of $\bar{A}^t(t-s+1,t)$ on $Y^t$ for $t \in \mathcal{T}_s$ based on MSMs for the $t$-specific full data associated with $O^t$. Figure 5 illustrates how the $t$-specific counterfactual framework can be viewed as a collection of conventional counterfactual sub-frameworks with distinct definition of the outcome, treatment and baseline covariates. These conventional counterfactual sub-frameworks differ from the conventional counterfactual framework in the sense that the treatment history is of size $s \neq K + 1$ and the outcome of interest is no longer time-dependent.

### 3.3. Three estimators of HRMSM parameters

We show in the appendix that an HRMSM corresponds to a $t$-specific MSM and the intersection of $t$-specific MSMs in a stratified and pooled analysis, respectively. This result implies that estimating functions for HRMSM parameters are $t$-specific MSM estimating functions (stratified analysis) or sums of $t$-specific MSM estimating functions (pooled analysis). Thus, we can easily extend the three estimators developed for MSM parameters to HRMSM parameters: the Inverse Probability of Treatment Weighted (16, 17, 19), the G-computation (12, 13, 3, 15, 29, 28, 8, 9) and Double Robust (18, 26, 7) estimators. To illustrate how the definition and implementation of MSM estimators are extended to HRMSM parameters based on the theoretical results in the appendix, we focus on IPTW estimation which has been favored in real-life applications.

### 3.4. IPTW estimator: definition and implementation

We present the two classes of IPTW estimators for HRMSM parameters in the stratified and pooled analyses based on a parametric HRMSM approach. The reader can easily extend these estimators to the nonparametric approach based on the MSM literature (10, 8).

#### 3.4.1. Stratified analysis

In this analysis, causal effects are modelled separately for each time point $t \in \mathcal{T}_s$, i.e., one separately investigates the causal effects on the outcomes of interest, $Y(t+1)$ for $t \in \mathcal{T}_s$, through the estimation of distinct causal parameters $\beta_t(F_X \mid \cdot)$ for $t \in \mathcal{T}_s$ defined based on $l = \text{Card}(\mathcal{T}_s)$ distinct models $m_t$ for $t \in \mathcal{T}_s$.

Under this model, the investigation of the causal effects of interest based on parametric models corresponds to the estimation of the $l$ distinct causal parameters $\beta_t^* \equiv \beta_t^*(F_X, g \mid m_t, g)$ defined such that:

$$E_{F_X,g}(Y_{\bar{a}(t-s+1,t)}(t+1) \mid V(t-s+1)) = m_t(\bar{a}(t-s+1,t), V(t-s+1) \mid \beta_t^*) \text{ for } t \in \mathcal{T}_s.$$

For a given $t$, the class of IPTW estimating functions for $\beta_t^*$ with nuisance parameter $g$ is defined as:

$$\left\{ D_h(O \mid g, \beta) = \frac{h(\bar{A}(t-s+1,t), V(t-s+1))\epsilon(\beta)}{\prod_{j=t-s+1}^{t} g(A(j) \mid \bar{A}(j-1), \bar{L}(j))} : h \in \mathcal{H} \right\},$$

where:

- $\varepsilon(\beta) = Y(t+1) - m_t(\bar{A}(t-s+1,t), V(t-s+1) \mid \beta)$,
- $\mathcal{H}$ is a set of non-null functions of $\bar{A}(t-s+1,t)$ and $V(t-s+1)$.

We denote estimators of $g$ and $h$ with $g_n$ and $h_n$.

We define the Experimental Treatment Assignment (ETA) assumption (7) for estimating $\beta_t^*$ as follows:

$$\max_{\bar{a}(t^-,t) \in \mathcal{A}_V(t^-,t)} \frac{h(\bar{a}(t^-,t), V(t^-))}{\prod_{j=t^-}^{t} g(a(j) \mid \bar{A}(t-s), \bar{a}(t^-, j-1), \bar{L}(j))} < \infty \ F_X - a.e,$$

where $t^- \equiv t-s+1$ and $\mathcal{A}_V(t-s+1, t)$ is the set of possible treatment regimens between time points $t-s+1$ and $t$, i.e. the support of the conditional distribution of $\bar{A}(t-s+1, t)$ given $\bar{A}(t-s)$ and $V$, $g(\bar{A}(t-s+1, t) \mid \bar{A}(t-s), V)$.

The IPTW estimator of $\beta_t^*$ is defined as the solution of the estimating equation associated with the observed data $O$ and the IPTW estimating function at $g_n$:

$$\sum_{i=1}^{n} D_{h_n}(o_i \mid g_n, \beta) = 0$$

Under regularity conditions, the IPTW estimator of $\beta_t^*$ is consistent and asymptotically linear if the ETA assumption holds and if $g_n$ is a consistent estimator of $g$.

In practice, the implementation of the IPTW estimator of the HRMSM parameter $\beta_t^*$ is identical to the implementation of MSM parameters in a stratified analysis with the difference that the treatment of interest for outcome $Y(t+1)$ is limited to $\bar{A}(t-s+1, t)$ instead of $\bar{A}(t)$. In other words, the IPTW estimate can be obtained with a weighted least squares regression of $Y(t+1)$ on $\bar{A}(t-s+1, t)$ and $V(t-s+1)$ based on the parametric model $m_t$ and weights inversely proportional to the estimated treatment mechanism. Like in MSM estimation, the numerator of the weights is defined by the choice for $h$ (e.g. stabilized versus unstabilized weights (19)). Unlike MSM estimation, note that the treatment mechanism is limited to treatments assigned between time points $t-s+1$ and $t$: $\bar{A}(t-s+1, t)$ and not all time points between baseline and $t$: $\bar{A}(t)$.

### 3.4.2. Pooled analysis

In this analysis, causal effects are modelled simultaneously for each time point $t \in \mathcal{T}_s$, i.e. change of the causal effect on the outcome over time is represented by a smooth function of time: one simultaneously investigates the causal effects on $Y(t+1)$ for $t \in \mathcal{T}_s$, through the estimation of a single causal parameter $\beta(F_X, g \mid \cdot)$ defined based on a single model $m(t, \bar{a}(t-s+1, t), V(t-s+1) \mid \beta)$.

Under this model, the investigation of the causal effects of interest based on parametric models corresponds to the estimation of the single causal parameter $\beta^* = \beta^*(F_X, g \mid m, g)$ defined such that

$$E_{F_X, g}(Y_{\bar{a}(t-s+1,t)}(t+1) \mid V(t-s+1)) = m(t, \bar{a}(t-s+1,t), V(t-s+1) \mid \beta^*)$$

for $t \in \mathcal{T}_s$.

For a given $t$, the class of IPTW estimating functions for $\beta^*$ with nuisance parameter $g$ is defined as:

$$\left\{ D_h(O \mid g, \beta) = \sum_{t \in \mathcal{T}_s} \frac{h(t, \bar{A}(t-s+1,t), V(t-s+1))\epsilon(\beta)}{\prod_{j=t-s+1}^{t} g(A(j) \mid \bar{A}(j-1), \bar{L}(j))} : h \in \mathcal{H} \right\},$$

where:

- $\varepsilon(\beta) = Y(t+1) - m(t, \bar{A}(t-s+1,t), V(t-s+1) \mid \beta)$,
- $\mathcal{H}$ is a set of non-null functions of $t$, $\bar{A}(t-s+1, t)$ and $V(t-s+1)$.

We denote estimators of $g$ and $h$ with $g_n$ and $h_n$.

We define the Experimental Treatment Assignment (ETA) assumption (7) for estimating $\beta^*$ as follows:

$$\max_{t \in \mathcal{T}_s} \max_{\bar{a}(t^-,t) \in \mathcal{A}_V(t^-,t)} \frac{h(t, \bar{a}(t^-,t), V(t^-))}{\displaystyle\prod_{j=t^-}^{t} g(a(j) \mid \bar{A}(t-s), \bar{a}(t^-,j-1), \bar{L}(j))} < \infty \; F_X - a.e,$$

where for each $t$, $t^- \equiv t-s+1$ and $\mathcal{A}_V(t-s+1, t)$ is the set of possible treatment regimens between time points $t-s+1$ and $t$, i.e. the support of the conditional distribution of $\bar{A}(t-s+1, t)$ given $\bar{A}(t-s)$ and $V$, $g(\bar{A}(t-s+1, t) \mid \bar{A}(t-s), V)$.

The IPTW estimator of $\beta^*$ is defined as the solution of the estimating equation associated with the observed data $O$ and the IPTW estimating function at $g_n$:

$$\sum_{i=1}^{n} D_{h_n}(o_i \mid g_n, \beta) = 0$$

Under regularity conditions, the IPTW estimator of $\beta^*$ is consistent and asymptotically linear if the ETA assumption holds and if $g_n$ is a consistent estimator of $g$.

In practice, the implementation of the IPTW estimator of the HRMSM parameter $\beta^*$ is identical to the implementation of MSM parameters in a pooled analysis with the difference that the treatment of interest for outcome $Y(t+1)$ is limited to $\bar{A}(t-s+1,t)$ instead of $\bar{A}(t)$. In other words, the IPTW estimate can be obtained with a pooled weighted least squares regression of $Y(t+1)$ for all $t \in \mathcal{T}_s$ on $\bar{A}(t-s+1,t)$ and $V(t-s+1)$ based on the parametric model $m$ and weights inversely proportional to the estimated treatment mechanism. Like in MSM estimation, the numerator of the weights is defined by the choice for $h$ (e.g. stabilized versus unstabilized weights (19)). Unlike MSM estimation, note that the treatment mechanism is limited to treatments assigned between time points $t-s+1$ and $t$: $\bar{A}(t-s+1,t)$ and not the entire treatments between baseline and time point $t$: $\bar{A}(t-s)$.

## 4. When and why prefer HRMSM-based versus MSM-based causal inference in practice?

### 4.1. MSM parameters: interpretation and causal effect representation

MSMs were introduced as a class of full data models which define parameters based on a feature of the marginal distribution of the following counterfactual outcomes: $Y_{\bar{a}(t)}(t+1)$ possibly conditional on the baseline covariates $V$. Typically and specifically in this article, one is interested in average causal effects per stratum $V$ of the population which can be represented by causal parameters defined by MSMs for $E_{F_X}(Y_{\bar{a}(t)}(t+1) \mid V)$ for $t \in \mathcal{T}_s$. We denote a causal parameter defined by an MSM with $\beta_t(F_X \mid \cdot)$ to indicate that it is a mapping from the space of full data distribution $F_X$ to the space of real numbers and that this mapping is a function of modeling assumptions represented by $\cdot$.

Two approaches to causal inference based on MSM have been proposed. They provide different representations of causal effects with distinct causal parameters. Initially, a parametric MSM approach to causal inference was developed (14) that relies on correct specification of a parametric MSM. Recently, a new approach based on nonparametric MSM was introduced (10, 8) that does not require to assume a correctly specified MSM and that generalizes the definition of causal parameters. This later approach is more realistic if one believes that correct specification of a parametric MSM is unlikely in practice.

In addition, both MSM approaches can be based on either a stratified or a pooled analysis, i.e. distinct models, $m_t(\bar{a}(t), V \mid \beta_t)$, or a single model, $m(t, \bar{a}(t), V \mid \beta)$, for $E_{F_X}(Y_{\bar{a}(t)}(t+1) \mid V)$ for $t \in \mathcal{T}$ (8).

Independently of the MSM approach chosen (nonparametric versus parametric and pooled versus stratified), MSM parameters represent the causal effects of the treatment histories, $\bar{A}(t)$, on the outcomes, $Y(t+1)$, for $t \in \mathcal{T}$. Note that this implies that in MSM-based causal inference, the causal effect of the treatment on the outcome collected at time point $t$ is always investigated for a treatment history of size $t$. As a result, the causal effect on the outcome collected at time $t$, $Y(t)$, is defined based on larger treatment histories as $t$ increases, i.e. as the outcome is collected later in the longitudinal study. This feature of this causal analysis is illustrated in figure 1. The figure illustrates how causal effects are investigated in practice based on an MSM approach for a study where $K = 5$, i.e. where the observed data is:

$$O = (L(0), A(0), L(1), A(1), L(2), A(2), L(3), A(3), L(4), A(4), L(5), A(5), L(6)).$$

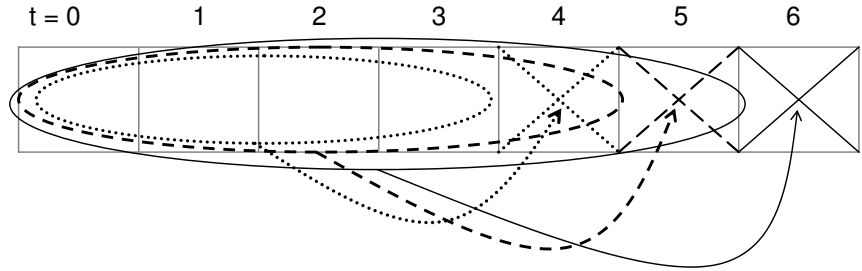

Fig 1: Illustration of the MSM representation of causal effects in a longitudinal study with time-dependent outcomes. Each ellipse indicates the relevant time-span over which the treatment's effect on a given outcome is investigated. Each cross indicates an outcome and each arrow represents one of the studied effects between treatment histories (ellipses) and outcomes (crosses).

For instance, the causal effect of the treatment on the outcome collected at time point $t = 6$ is investigated for a treatment history of size 6, $\bar{A}(5)$, as represented on figure 1 by the ellipse that covers time point 0 to 5 and the arrow that connects the ellipse to the outcome collected at time point 6 represented by a cross.

The investigation of causal effects with MSMs (see figure 1) raises three *potential* concerns that are likely to become more significant for longitudinal studies with longer follow-ups:

- computational intractability when proceeding to MSM estimation
- a disbelief about the subject-matter relevance of the causal effects investigated
- a statistical power problem

The first issue is best illustrated with the implementation of the G-computation estimator and was fully developed in previous work (8, 9). More generally, implementation of MSM estimators, like the IPTW estimator, becomes less practicable as the follow-up time increases.

The second issue can easily be illustrated with the following hypothetical study. Consider a longitudinal study during which individuals are treated or not every day over three months (90 days) with a new medication for headache relief and monitored for headache symptoms. Now consider the last outcome, $Y(90)$, collected after a treatment history of 89 days, $\bar{A}(89)$. In MSM-based causal inference, the investigation of treatment effect on the last outcome measured would be based on the estimation of a causal parameter representing the effect of the treatment history $\bar{A}(89)$ on $Y(90)$. Most investigators would argue that looking a the effect of such a medication taken 3 months prior outcome report is likely to be of little interest since: 1) the relief effect of such medication usually does not carry over such a long period of time; and 2) because the drug effect that is pursued, i.e. of interest, for such a treatment is a short-term relief. In other words, investigation of the effect of a headache reliever absorbed 3 months prior outcome report is not of primary interest. It may also often make no sense to study such a long lag effect in practice if the treatment is known to act over a short-term time scale exclusively. Note that a natural way to overcome this problem in practice would be to investigate short-term effects based on a parametric MSMs only involving the latest treatments received before the outcome is collected. For example a parametric MSM for $Y_{\bar{a}(89)}(90)$ may only rely on the last two treatments absorbed to explain the outcome, e.g. $E(Y_{\bar{a}(89)}(90)) = \beta_0 + \beta_1 a(88) + \beta_2 a(89)$. Note however that such an MSM assumes that treatment before time 88 has no effect on the treatment which may be incorrect and estimation of $\beta_1$ and $\beta_2$ should still rely on the complete treatment history between time point 0 and 89. We argue in this article that the MSM approach can be improved to better identify causal effects that are truly of interest from a subject-matter point of view without making restrictive assumptions often incorrect in practice. The HRMSM proposed in this article addresses this issue and mitigates the other two concerns discussed here.

The third issue is not illustrated in this article with a concrete example but has been discussed in previous work (2, 6) and is related to the Experimental Treatment Assignment assumption (7). In short, in statistical analyses based on MSMs, the longer the treatment history, i.e. the study follow-up, the more complex the description of the treatment effect becomes. Thus, it is likely that the information required to understand long-term effects will be very important and beyond the reach of most investigators. Concretely, the decrease in statistical power can be explained by an increasing practical violation of the ETA assumption as the treatment history increases. As a result, even when it is sensible to investigate causal effects based on MSMs in longitudinal studies with long follow-up, investigators may still wish to revise their study aims and lower their research ambitions for the sake of practicability. Note that it is this statistical power issue which initially motivated the introduction of HRMSMs with failure-time data (2, 6). The authors indeed noted that MSM estimation with the IPTW estimator suffered from a lack of precision for long treatment histories due to an increasing variability in the weights. This phenomemon can be explained by practical violation of the ETA assumption when the treatment history is long (7).

### *4.2. HRMSM parameters: interpretation and causal effect representation*

Both the parametric and nonparametric MSM (14, 10, 8) approaches that have been proposed in causal inference can be directly extended to HRMSM-based approaches. Similarly, the corresponding parametric and nonparametric HRMSM approaches to causal inference provide different representations of causal effects with distinct causal parameters.

In addition, both HRMSM approaches can be based on either a stratified or a pooled analysis, i.e. distinct models, $m_t(\bar{a}(t-s+1,t), V(t-s+1) \mid \beta_t)$, or a single model, $m(t, \bar{a}(t-s+1,t), V(t-s+1) \mid \beta)$, for $E_{F_X,g}(Y_{\bar{a}(t-s+1,t)} \mid V(t-s+1))$ for $t \in \mathcal{T}_s$ (8).

Independently of the HRMSM approach chosen (nonparametric versus parametric and pooled versus stratified), the causal effect of the treatment on the outcome is always investigated for a treatment history of fixed size, $s$, in HRMSM-based causal inference. As a result, the causal effect on the outcome collected at time $t$, $Y(t)$, is defined based on a fixed treatment history even as $t$ increases, i.e. as the outcome is collected later in the longitudinal study. This feature of this causal analysis is illustrated in figure 2. It illustrates how causal effects are investigated in practice with HRMSMs for a study where $K = 5$ and $s = 2$, i.e. where the data collected is:

$$(L(0), A(0), L(1), A(1), L(2), A(2), L(3), A(3), L(4), A(4), L(5), A(5), L(6)).$$

For instance, the causal effect of the treatment on the outcome collected at time point $t = 6$ is investigated for a treatment history of size 2: $(A(4), A(5))$ as represented on figure 2 by the ellipse that covers time points 4 to 5 and the arrow that connects the ellipse to the outcome collected at time point 6 represented by a cross.

The following three practical considerations may often lead investigators to prefer a causal analysis based on HRMSMs above an analysis based on MSMs for causal inference problems with longitudinal data:

- Computational tractability,
- Causal effect representation that is most relevant to public health research,
- Statistical power.

Support for these claims can be found in section 1 and 4.1 where the limitations of MSM-based causal inference are underscored.

## 5. Illustration with an air pollution study

### *5.1. Specific aims and data*

We illustrate the application of HRMSMs with one of the air pollution study already described in section 1. Its primary goal is to investigate the extent to which reductions in ambient air pollution consequent to regulations propagated since 1980 by the California Air Resources Board to reduce air pollution in the Los Angeles (LA) Basin lead to measurable health benefits. In this example, we will specifically focus on the reduction of ozone levels and its effect on health as measured by asthma-related hospital discharges in children ages birth to 19 years[1].

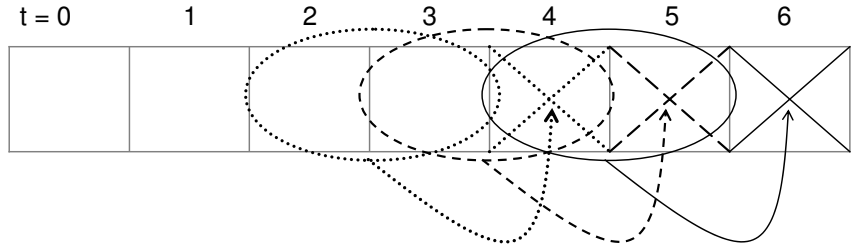

Fig 2: Illustration of the HRMSM representation of causal effects in a longitudinal study with time-dependent outcomes and $s = 2$. Each ellipse indicates the relevant time-span over which the treatment's effect on a given outcome is investigated. Each cross indicates an outcome and each arrow represents one of the studied effects between treatment histories (ellipses) and outcomes (crosses).

The experimental units are $10 \times 10$ km geographical areas of the LA Basin. Data were assembled on a total of $n = 195$ such grids over 72 ($K = 70$) quarters to answer the aforementioned question of interest. For each grid, $i$, and each quarter $t$, the data set contains measurements for 1) the exposure of interest: ozone levels denoted with $A(t)$, 2) the outcome of interest: asthma-related hospital discharge counts denoted with $C(t)$ and the total number of hospital discharges $N(t)$ in children ages birth to 19 years, and 3) 56 socioeconomic and demographic covariates denoted with $L(t)$ (e.g. age, income, racial and gender structures as measured by ratios of the grid census).

### 5.2. Model and parameter of interest

To address the question of interest, we believe that an MSM-based statistical strategy to analyze these data is most natural. Indeed, the overall subject-matter objective is two-fold: 1) we wish to verify that exposure to ozone indeed lead to averse public health impact, i.e. we wish to evaluate the causal effect of ozone on the asthma-related hospital discharge rate; and 2) after verification of the averse nature of this effect, we wish to evaluate the public health benefits related to (caused by) decrease in ozone levels over 20 years (see figure 3), i.e. we wish

---

[1]Results are for heuristic purposes only and are not presented as a definitive analysis on the subject.

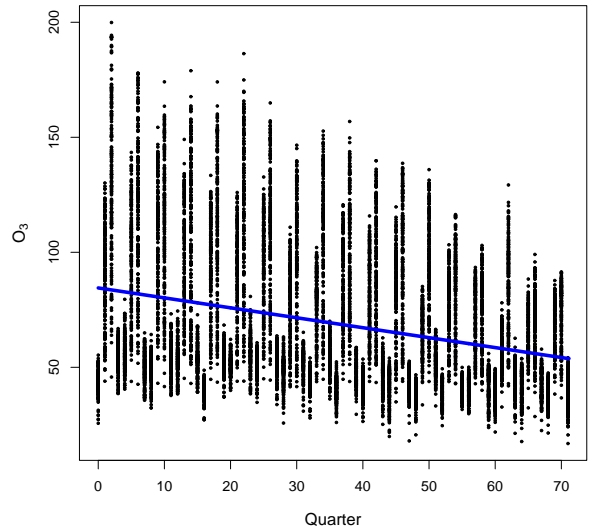

Fig 3: Plot representing the decrease of ozone levels over 72 quarters for the 195 grid units. The linear curve corresponds with the fit of a linear regression model of ozone levels against time.

to contrast what would the population outcomes be had the ozone levels not decreased over time.

A central issue that emerges is the relevant exposure time for investigation of the exposure effect on the outcome rate: it does not seem reasonable to extend the exposure period much beyond the 12 months prior to a given quarter. We indeed believe that most of the effect of ozone can be captured within this time frame. This choice of a time frame is reinforced by the fact that the experimental units are not individuals but geographical areas: the population constituting each grid is constantly changing and major modifications in the population structures occurred over the 20-year period studied. Our focus is, thus, on the sub-acute effects: we consider the effect of one-year exposure to ozone during which time the population in each grid is relatively constant. Even, if the effects of ozone were to extend much beyond this time frame, we believe that we would lack power to identify such effects with our finite data ($n = 195$). In addition, software development is too limited to date to engage in a more complex analysis where we would investigate the ozone levels over large histories. Consequently, an HRMSM is an appealing statistical tool for the analyses of these data.

We propose to use a pooled, binomial logistic HRMSM with $s = 4$ (i.e., we consider the effect of 4 quarters of exposure to ozone) to represent the effect of interest in this analysis, i.e. we assume that the conditional distribution of $C_{\bar{a}(t-3,t)}(t+1)$ conditional on $N_{\bar{a}(t-3,t)}(t+1)$ follows a binomial distribution:

$C_{\bar{a}(t-3,t)}(t+1) \sim \mathcal{B}(N_{\bar{a}(t-3,t)}(t+1), p_{\bar{a}(t-3,t)}(t+1))$ with:

$$E\left(\frac{C_{\bar{a}(t-3,t)}(t+1)}{N_{\bar{a}(t-3,t)}(t+1)} \mid N_{\bar{a}(t-3,t)}(t+1)\right)$$

$$= \text{logit}\Bigg(\beta_0 + \beta_1 \text{mean}(\bar{a}(t-3,t)) + \beta_2 f_1(t) + \beta_3 f_2(t)$$

$$+ \beta_4 \text{mean}(\bar{a}(t-3,t)) f_1(t) + \beta_5 \text{mean}(\bar{a}(t-3,t)) f_2(t)$$

$$+ \beta_6 \text{mean}(\bar{a}(t-3,t)) f_1(t) f_2(t)\Bigg),$$

where:

- $t \in \mathcal{T}_4 = \{3, \ldots, 70\}$
- $\text{logit}(x) = \frac{1}{1+\exp(-(x))}$,
- $\text{mean}(\bar{a}(t-3,t)) = \frac{a(t-3)+a(t-2)+a(t-1)+a(t)}{4}$,
- $f_1(t)$ is a mapping from $t$ to the year number associated with quarter $t+1$, and
- $f_2(t)$ is a mapping from $t$ to the season associated with quarter $t+1$ (two 6-month seasons are considered in the analysis only: winter and summer).

This model implies the following:

$$\log\left(\frac{p_{\bar{a}(t-3,t)}(t+1)}{1-p_{\bar{a}(t-3,t)}(t+1)}\right) = \beta_0 + \beta_1 \text{mean}(\bar{a}(t-3,t)) + \beta_2 f_1(t) + \beta_3 f_2(t)$$

$$+ \beta_4 \text{mean}(\bar{a}(t-3,t)) f_1(t) + \beta_5 \text{mean}(\bar{a}(t-3,t)) f_2(t)$$

$$+ \beta_6 \text{mean}(\bar{a}(t-3,t)) f_1(t) f_2(t),$$

because $E\left(\frac{C_{\bar{a}(t-3,t)}(t+1)}{N_{\bar{a}(t-3,t)}(t+1)} \mid N_{\bar{a}(t-3,t)}(t+1)\right) = p_{\bar{a}(t-3,t)}(t+1)$.

### 5.3. Estimation

To fit the HRMSM, i.e. estimate $\beta = (\beta_1, \ldots, \beta_6)$, we chose the G-computation estimator. A detailed algorithm for the implementation of the G-computation estimator has been described in previous work (8, 9). This implementation relies, in practice, on a data reduction step based on a generalized linear model for the treatment mechanism. Based on the reduced data set, the G-computation estimate of $\beta$ is obtained through Monte Carlo simulation based on a generalized linear model for the $Q_{F_X}$ part of the likelihood. All models involved in this estimation procedure were obtained through model selection. We utilized a recently developed data-adaptive model selection procedure (23, 25) based on cross-validation and an aggressive model search algorithm known as the D/S/A algorithm (21).

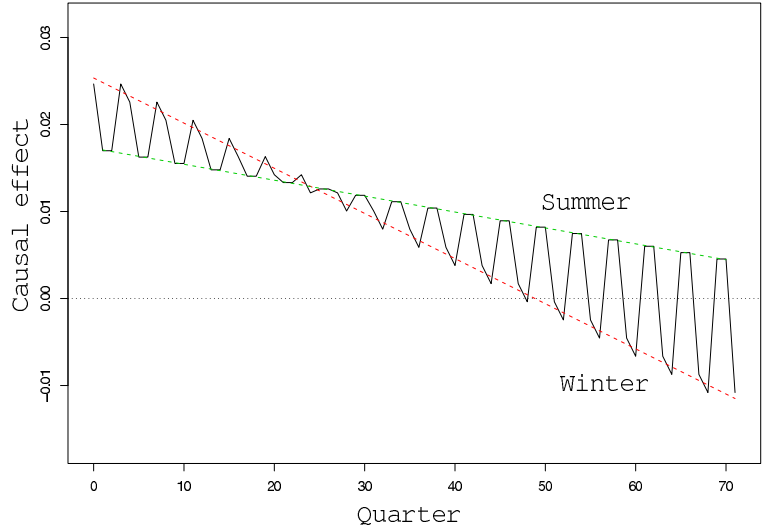

Fig 4: Effect of one-year exposure to ozone on asthma-related hospital discharges over time (curvilinear line). The straight lines represent the season-specific causal effects. Each point above the $y = 0$ line is interpreted as a decrease in the odds of asthma-related hospital discharge incidence caused by a decrease of one unit in the average ozone quarterly levels over one year.

### *5.4. Results and interpretation*

Note that the results presented are preliminary[2]. Figure 4 represents a summary of the results with a plot of the temporal change of the effect of one-year exposure to ozone on the asthma-related hospital discharge rate as measured by $\hat{\beta}_1 + \hat{\beta}_4 f_1(t) + \hat{\beta}_5 f_2(t) + \hat{\beta}_6 f_1(t) f_2(t)$ where $\hat{\beta}$ represents the G-computation point estimate of $\beta$. Each point above the $y = 0$ line on figure 4 is interpreted as a decrease in the odds of asthma-related hospital discharge incidence caused by a decrease of one unit in the average ozone quarterly levels over one year. In addition, the unit effect of ozone on asthma-related hospital discharge incidence decreases linearly over time. Note also that a portion of the winter section of the curve is below the $y = 0$ line which indicates a decrease in asthma-related admissions per unit increase of ozone levels. This decrease could represent uncontrolled confounding, measurement errors and/or model misspecification. We did not explore this finding, since results have been presented to illustrate the method. Also, we do not provide confidence intervals which can be obtained with the bootstrap as is appropriate with the MSM methodology.

---

[2]Results are for heuristic purposes only and are not presented as a definitive analysis on the subject.

## 6. Discussion

In this paper, we have laid out the formal statistical framework for a new class of MSMs, HRMSMs. These models were initially introduced (2, 6) as alternative causal inference tools to MSMs in order to improve statistical power when investigating the causal effects of a treatment history on a health outcome. We developed an extension of the conventional counterfactual framework that we called the *t*-specific counterfactual framework. This framework was introduced solely as a statistical artifice to provide the rigorous mathematical framework to develop consistent estimators of HRMSM parameters with minimal effort: the IPTW, G-computation and DR estimators. We have shown that these estimators of HRMSM parameters are consistent under the same model assumptions commonly adopted in the conventional counterfactual framework: existence of counterfactuals, consistency, time-ordering and sequential randomization assumptions.

In addition, we further argued in this article, based on practical considerations (computational tractability, relevance of the effect investigated, statistical power), that an HRMSM-based causal inference strategy may often be better suited for public health research than an MSM approach where the effect of an exposure history is modeled as a function of exposures experienced only in the last $s$ periods prior outcome collection. Such an MSM approach implicitely relies on the assumption that exposures experienced before the last $s$ periods had no effect on the outcome. This assumption typically does not hold in most applications. The HRMSM approach makes no such assumption, and so is more circumspect in its implications. We believe these considerations should motivate the application of this methodology in many epidemiological and clinical studies. We now discuss the decision making about the history size, $s$, in real-life applications of the proposed HRMSM-based causal analysis.

Decision about the value for the history size, $s$, should be based on the combination of considerations about the analysis aims and *a priori* knowledge about the problem being studied as illustrated in section 5. This decision however cannot ignore practical considerations such as implementation issues and statistical power concerns.

For instance, if a longitudinal study aims to investigate the causal effect of a new medication for headache relief whose action is likely not to carry over time beyond a few hours, then it will not make sense to choose a history size that extends well beyond the known lag effect of similar medication. Even when the treatment effect of interest is likely to carry over long periods of time, the subject-matter focus may be the investigation of short-term effects, in which case the investigators should consider small values for $s$. In addition, note that the larger $s$, the more complex the causal effect of interest that is captured in the time interval represented by $s$. As a result, the statistical power to investigate the causal effect of interest will likely decrease when considering larger history sizes $s$. Moreover, the larger $s$ is, the less computationally tractable will HRMSM estimation be.

Nevertheless, choosing a suitable value for $s$ based on these guidelines remain subjective and may lead to two scenarios: 1) the chosen history size, $s$, is larger than the maximum time interval over which the treatment of interest has an effect on the outcome; and 2) the chosen history size, $s$, is smaller than the maximum time interval over which the treatment of interest has an effect on the outcome. In the first case scenario, the model selection procedure for MSMs proposed by van der Laan and Dudoit (2003) (24) can be used to identify the smaller component of the treatment history that is causally relevant. The second case scenario is most likely to occur in practice, since statistical power and implementation considerations often will prevent investigation of the causal effects of treatment histories that are too long. An HRMSM-based causal analysis still will provide valuable answers to the public health questions of interest, based on the available data even if the causal effect of the treatment over time will not be completely described (e.g. the maximum lag effect will remain unknown).

We would also like to note that a variant of HRMSMs can be used for possibly discrete failure-time outcomes. For a discrete failure-time outcome, this would be a model for $E_{F_X,g}(Y_{\bar{a}(t-s+1,t)}(t+1) \mid Y_{\bar{a}(t-s+1,t)}(t), V(t-s+1))$; this is a model for the probability of failure at time $t+1$ under a given restricted treatment regimen conditional on having survived through time $t$ under that regimen and possibly $t$-specific baseline covariates. Such models pose two problems for causal interpretation. First, as is generally true with comparisons of hazards, even in randomized trials, the comparison of two treatment regimens' effects at a given time point is not based on comparing outcomes for a common set of individuals; rather, it is based on comparing failure among survivors under one treatment regimen to failure for another set of survivors under the other treatment regimen. If subjects susceptible to harmful effects of a given treatment regimen fail early after being treated with that regimen, the depletion of susceptibles may lead to hazards under that treatment regimen being the same or even lower than in its absence. In simple randomized trials or with MSMs, this problem with interpretation can be mitigated by using the hazards under a given treatment regimen to compute and compare treatment-specific survival functions instead of treatment specific hazard functions. The comparisons of treatment-specific survival functions are comparisons for a common group of individuals under different treatment regimens, starting from the beginning of the regimen. In this variant of HRMSMs, treatment-specific survival curves cannot be computed, because unlike MSMs, hazards at different time points $t$ in HRMSMs do not correspond with a single treatment-specific hazard of interest and can thus not be mapped into a treatment specific survival curve. This is particularly obvious if the HRMSMs condition on covariates $V(t-s+1)$ since the set of covariates $V(t-s+1)$ will not be identical over time. These concern do not apply to HRMSMs for non-survival outcomes.

Finally, we would like to clarify the nomenclature for models that have recently been proposed by different authors to essentially designate the same class of models: partially MSMs (6), History-Adjusted MSMs (27), and History-Restricted MSMs. The difference in the nomenclature reflects the heterogeneity in the subject-matter problems for which different authors developed a modeling approach based on the same class of MSMs. Unlike partially MSMs and HRMSMs, the primary motivation for the development of History-Adjusted MSMs was the inability to incorporate modification of the causal effect of treatment by time-varying covariates in MSMs for clinical decision making. It is interesting to note that solutions to these different problems can be found in the same class of MSMs. This confluence should further underscore the importance of the proposed modeling approach in practice.

# APPENDIX

## Appendix A: HRMSM estimation: the time-specific counterfactual framework

In this section, we introduce the time-specific ($t$-specific) counterfactual framework which can be viewed as an extension of the conventional counterfactual framework on which MSM-based causal inference is based (see sections 2.1 and 2.2). This latter mathematical construct provided the rigorous framework to define, identify and estimate MSM parameters with the full and observed data based on a sufficient set of assumptions developed in section 2.2. We introduce the $t$-specific counterfactual framework because it allows us to generalize the MSM estimation procedures to HRMSM estimation procedures with minimum effort. In addition, we present in this section the sufficient assumptions for estimation of HRMSM parameters in the $t$-specific counterfactual framework.

### A.1. Data structures

In this section, we adopt the notations introduced in the previous section to represent the treatments, covariates and outcomes collected at each time point $t = 0, \ldots, K+1$ on each of the $n$ experimental units: $A(t), L(t), Y(t)$ respectively. We also adopt the notation $V(t-s+1)$ to designate a subset of $(\bar{A}(t-s), \bar{L}(t-s+1))$.

In the $t$-specific counterfactual framework, the representation of the data collected during a longitudinal study between time points 0 and $K+1$ is based on a user-specified choice of a fixed treatment history size $s > 0$. We already discussed the interpretation of this parameter, $s$. In section 6, we discuss the decision making about its value in practice.

In the $t$-specific counterfactual framework and for a given treatment history size $s$, the data are represented as $n$ i.i.d realizations of $K - s + 2$ data structures:

$$
\begin{aligned}
O^t \;=\; & (L^t(t-s+1), A^t(t-s+1), L^t(t-s+2), A^t(t-s+2), \ldots, L^t(t), \\
& \hspace{6cm} A^t(t), L^t(t+1)) \\
\;=\; & (\bar{A}^t(t-s+1,t), \bar{L}^t(t-s+1,t+1)) \sim P^t,
\end{aligned}
$$

for all $t$ such that $s - 1 \leq t \leq K$ and where

- $P^t$ represents the distribution of the stochastic process $O^t$ referred to as one of the *t*-specific observed data,
- $\bar{A}^t(t - s + 1, t)$ represents the *t*-specific treatment process defined as

$$A^t(j) \equiv A(j) \text{ for all } j \text{ such that } t - s + 1 \leq j \leq t, \quad (1)$$

- $\bar{L}^t(t - s + 1, t + 1)$ represents the *t*-specific covariate process defined by:

  a)  $L^t(j) \equiv L(j)$ for all $j$ such that $t - s + 1 < j \leq t + 1$      (2)

  b)  $L^t(t - s + 1) \equiv (\bar{A}(0, t - s), \bar{L}(0, t - s + 1))$      (3)

  In other words, we have $\bar{L}^t(t - s + 1, t + 1) \equiv (\bar{A}(t - s), \bar{L}(t + 1))$.

Furthermore, we define $V^t$ as a subset of the baseline covariates in the *t*-specific observed data, $O^t$, $V^t \subset L^t(t - s + 1)$, such that:

$$V^t \equiv V(t - s + 1). \quad (4)$$

We define $Y^t$ as the *t*-specific outcome of interest, $Y^t \in L^t(t + 1)$, such that:

$$Y^t \equiv Y(t + 1). \quad (5)$$

In addition, we denote with $\mathcal{T}_s$ the set of time points $t$ such that the outcome $Y^t$ is of interest. We have $\mathcal{T}_s \subset \{s - 1, \ldots, K\}$. Typically we will have $\mathcal{T}_s = \{s - 1, \ldots, K\}$. Note that the outcome of interest maybe time dependent, e.g. $\mathcal{T}_s = \{s - 1, \ldots, K\}$, while the *single* outcome of interest in the *t*-specific observed data is $Y^t \equiv Y(t + 1)$, i.e. $Y(t)$ is not an outcome of interest in $O^t$ but only in $O^{t-1}$.

Like in the conventional counterfactual framework, the question of interest is to investigate the causal effect of treatment $A$ on the time-dependent outcome, $Y \in L$. In the *t*-specific counterfactual framework, this problem is addressed through the investigation of the causal effects of the treatment histories $\bar{A}^t(t - s + 1, t)$ on the outcomes $Y^t \in L^t(t + 1)$ for all $t \in \mathcal{T}_s$.

We want to underscore again that in this approach and for a given $t \in \mathcal{T}_s$, the outcome $Y^t$ is not defined as a time-dependent variable in the sense that it corresponds with a variable measured at a unique time-point, specifically the last time-point $t + 1$ associated with the corresponding *t*-specific observed data, $O^t$. That is why, although $Y^t \equiv Y(t + 1)$, we adopt a separate notation $Y^t$ to designate the outcome. It is not to be confused with the notation $Y$ introduced for the conventional counterfactual framework and which designates a time-dependent variable. Indeed, $Y(j)$ can be regarded in the *t*-specific counterfactual framework both as any covariate $Y(j) \in L^t(j)$ and the outcome $Y^t$ when $j = t + 1$. Similarly, note that we adopt a distinct notation, $A^t$, to unambiguously represent the treatment of interest in the *t*-specific observed data $O(t)$. This notation is not to be confused with $A$ which refer to a variable that can be regarded both as a covariate $L^t$ and a treatment variable $A^t$ in the $t$-specific counterfactual framework.

Note that if for example $s = 2$, then $Y(1)$ cannot be of interest since at time point 1 each unit has been treated with a treatment history of size 1 and it is thus not possible to look at the effect of a treatment history of size 2 on $Y(1)$. That is why we use the notation $\mathcal{T}_s$ to indicate that the set of outcomes of interest depends on the investigator's choice for $s$.

### A.2. Assumptions

In the $t$-specific counterfactual framework, we adopt the same set of assumptions as described for the conventional counterfactual framework with the exception that each assumption is made relative to each $t$-specific observed data of interest, $O^t$ for $t \in \mathcal{T}_s$. In other words, we make the following $t$-specific assumptions **for all $\mathbf{t} \in \mathcal{T}_\mathbf{s}$.**

**Existence of counterfactuals:** we assume the existence of the following $t$-specific treatment specific processes, $\bar{L}^t_{\bar{a}^t(t-s+1,t)}(t-s+1, t+1)$, also referred to as $t$-specific counterfactual processes, for every treatment regimen $\bar{a}^t(t-s+1, t) = (a(t-s+1), \ldots, a(t)) \in \mathcal{A}_V(t-s+1, t)$ where $\mathcal{A}_V(t-s+1, t)$ designates all possible treatment regimens between time points $t-s+1$ and $t$, i.e. the support of the conditional distribution of $\bar{A}(t-s+1, t)$ given $\bar{A}(t-s)$ and $V$, $g(\bar{A}(t-s+1, t) \mid \bar{A}(t-s), V)$. We denote the so-called $t$-specific full data process associated with $O^t$ with $X^t = (V, \bar{L}^t_{\bar{a}^t(t-s+1,t)}(t-s+1, t+1))_{\bar{a}^t(t-s+1,t) \in \mathcal{A}_V(t-s+1,t)}$ and its distribution with $F_{X^t}$.

Note that the existence of the $t$-specific counterfactual process $\bar{L}^t_{\bar{a}^t(t-s+1,t)}(t-s+1, t+1)$ for every treatment regimen $\bar{a}^t(t-s+1, t) \in \mathcal{A}_V(t-s+1, t)$ implies the existence of the $t$-specific counterfactual processes $\bar{L}^t_{\bar{a}^t(t-s+1,j)}(t-s+1, j+1) \equiv \bar{L}^t_{\bar{a}^t(t-s+1,j), A^t(j+1), \ldots, A^t(t)}(t-s+1, j+1) \subset X^t$ for every $j = t-s+1, \ldots, t-1$ and every treatment regimen $\bar{a}^t(t-s+1, j) = (a(t-s+1), \ldots, a(j)) \in \mathcal{A}_V(t-s+1, j)$ where $\mathcal{A}_V(t-s+1, j)$ designates all possible treatment regimens between time points $t-s+1$ and $j$, i.e. the support of the conditional distribution of $\bar{A}(t-s+1, j)$ given $\bar{A}(t-s)$ and $V$, $g(\bar{A}(t-s+1, j) \mid \bar{A}(t-s), V)$. We have $\mathcal{A}_V(t-s+1, j) = \{\bar{a}(t-s+1, j) : \exists\, \bar{a}'(t-s+1, t) \in \mathcal{A}_V(t-s+1, t)\ \bar{a}(t-s+1, j) = \bar{a}'(t-s+1, j)\}$ for $j = t-s+1, \ldots, t-1$ and $\mathcal{A}_V(t-s+1, j)$ is thus entirely defined by $\mathcal{A}_V(t-s+1, t)$. Similarly, $\mathcal{A}_V(t-s+1, t)$ is entirely defined by $\mathcal{A}_V(K)$.

**Consistency assumption:** at any time point $j$ such that $t-s+1 \leq j \leq t+1$, we assume the following link between the $t$-specific observed data and the $t$-specific counterfactuals: $L^t(j) = L^t_{\bar{A}^t(t-s+1,t)}(j)$. Under this assumption, we have: $O^t = (\bar{A}^t(t-s+1, t), \bar{L}^t_{\bar{A}^t(t-s+1,t)}(t-s+1, t+1)) \equiv \phi^t(\bar{A}^t(t-s+1, t), X^t)$, where $\phi^t$ is a specified function of the $t$-specific full data process $X^t$. This notation indicates that the problem can be treated as multiple (for each $t \in \mathcal{T}_s$) missing data problems. Indeed, for each $t \in \mathcal{T}_s$, only the $t$-specific counterfactual associated with the observed treatment $\bar{A}^t(t-s+1,t)$ is observed; the others are missing.

**Temporal Ordering assumption:** at any time point $j$ such that $t-s+1 \leq j \leq t+1$, we assume that any treatment-specific variable can only be affected by past treatments: $L^t_{\bar{a}^t(t-s+1,t)}(j) = L^t_{\bar{a}^t(t-s+1,j-1)}(j)$ for $j = t-s+1, \ldots, t+1$, where $L^t_{\bar{a}^t(t-s+1,t-s)}(t-s+1) = L^t(t-s+1)$. This assumption is typically implied by the data collection procedure: the covariate $L^t(t)$ is measured after $A^t(t-1)$ and before $A^t(t)$.

**Sequential Randomization Assumption (SRA):** at any time point $j$ such that $t-s+1 \leq j \leq t+1$, we assume that the $t$-specific observed treatment is independent of the $t$-specific full data given the $t$-specific data observed up to time point $j$: $A^t(j) \perp X^t \mid \bar{A}^t(t-s+1, j-1), \bar{L}^t(t-s+1, j)$. Under the SRA, the $t$-specific treatment mechanism, i.e. the conditional density or probability of $\bar{A}^t(t-s+1,t)$ given $X^t$, $g(\bar{A}^t(t-s+1,t) \mid X^t)$, becomes:

$$
\begin{aligned}
g(\bar{A}^t(t^-,t) \mid X^t) &= \prod_{j=t^-}^{t} g(A^t(j) \mid \bar{A}^t(t^-,j-1), X^t) \\
&\stackrel{SRA}{=} \prod_{j=t^-}^{t} g(A^t(j) \mid \bar{A}^t(t^-,j-1), \bar{L}^t(t^-,j)),
\end{aligned}
$$

where $t^- \equiv t-s+1$. The SRA implies coarsening at random (4) and thus the $t$-specific likelihood of the $t$-specific observed data factorizes into two parts: a so-called $F_{X^t}$ and $g^t$ part. The $F_{X^t}$ part of the likelihood only depends on the $t$-specific full data process distribution, and the $g^t$ part of the likelihood only depends on the $t$-specific treatment mechanism. As a consequence of this factorization of the $t$-specific likelihood under the SRA, we now denote the distribution of the $t$-specific observed data with $P_{F_{X^t}, g^t}$ and the likelihood of $O^t$ is:

$$
\mathcal{L}(O^t) \stackrel{SRA}{=} \underbrace{f(L^t(t-s+1)) \underbrace{\prod_{j=t-s+2}^{t+1} f(L^t(j) \mid \bar{L}^t(t-s+1,j-1), \bar{A}^t(t-s+1,j-1))}_{Q_{F_{X^t}}}}_{F_{X^t} \text{ part}} \underbrace{g(\bar{A}^t(t-s+1,t) \mid X^t)}_{g^t \text{ part}}.
$$

In addition, we denote the set of conditional densities or probabilities that define the $F_{X^t}$ part of the likelihood, except for $f(L^t(t-s+1))$ with $Q_{F_{X^t}}$.

### A.3. Equivalence between MSM parameters in the t-specific counterfactual framework and HRMSM parameters in the conventional counterfactual framework

Like in the conventional counterfactual framework, causal effects can be represented based on parameters defined by MSMs in the $t$-specific counterfactual framework. Indeed MSM approach can be applied to all $t$-specific observed data, $O^t$. We refer to an MSM associated with a given $t$-specific full data as a $t$-specific MSM. These $t$-specific MSMs are $t$-specific full data models, i.e. model of $F_{X^t}$, which define parameters based on a feature of the distribution of the following counterfactual outcomes: $Y^t_{\bar{a}^t(t-s+1,t)}$. Typically and specifically in this article, one is interested in average causal effects per stratum $V^t$ of the population which can be represented by causal parameters defined by MSMs for $E_{F_{X^t}}(Y^t_{\bar{a}^t(t-s+1,t)} \mid V^t)$ for $t \in \mathcal{T}_s$. We denote a causal parameter defined by such an MSM with $\beta_t(F_{X^t} \mid \cdot)$ to indicate that it is a mapping from the space of $t$-specific full data distribution $F_{X^t}$ to the space of real numbers and that this mapping is a function of modeling assumptions represented by $\cdot$.

We have by definition from (4): $V^t \equiv V(t-s+1)$ and we can show as follows that $Y^t_{\bar{a}^t(t-s+1,t)} = Y_{\bar{a}(t-s+1,t)}(t+1)$ for $t \in \mathcal{T}_s$:

$$
\begin{aligned}
Y^t_{\bar{a}^t(t-s+1,t)} &= Y^t_{\bar{a}(t-s+1,t)} \text{ from (1)} \\
&= Y_{\bar{a}(t-s+1,t)} \text{ from (5)}
\end{aligned}
$$

Thus we have $E_{F_{X^t}}(Y^t_{\bar{a}^t(t-s+1,t)} \mid V^t) = E_{F_X,g}(Y_{\bar{a}(t-s+1,t)}(t+1) \mid V(t-s+1))$ and $F_{X^t} = \psi(F_X, g)$ for some specified function $\psi$.

In general, one can show that HRMSM parameters defined in the conventional counterfactual framework, $\beta_t(F_X, g \mid \cdot)$ (see section 2.3), corresponds to MSM parameters defined in the $t$-specific counterfactual framework, $\beta_t(F_{X^t} \mid \cdot)$:

$$
\beta_t(F_X, g \mid \cdot) = \beta_t(F_{X^t} \mid \cdot) \tag{6}
$$

In addition, note that MSM parameters defined in the $t$-specific counterfactual framework are typically different from MSM parameters defined in the conventional counterfactual framework. Indeed the former are defined as a function of both $F_X$ and $g$ while the latter are defined as a function of $F_X$ only.

### A.4. Link between the conventional and t-specific counterfactual frameworks

Figure 5 illustrates based on an example of a longitudinal study with short follow-up the link between the longitudinal data representation in the conventional counterfactual framework and its representation in the time-specific counterfactual framework. Note that in the conventional counterfactual framework the data are approached as a single entity, $O$, in the sense that the treatment is defined once and for all as a history $\bar{A}(K)$ and the outcome is time-dependent, $Y(t) \in L(t)$. On the other hand, in the $t$-specific counterfactual framework the data are viewed as layers of separate entities, $O^t$. For each $O^t$, the treatment and outcome of interest are redefined along with the baseline covariates (highlighted in yellow/light gray on figure 5). Note that for each $O^t$, the outcome is no longer time-dependent but correspond with the last outcome collected, $Y^t \in L^t(t+1)$, (highlighted in orange/dark gray on figure 5) and the treatment history size is fixed to a user-specified value $s = 2$. In the conventional counterfactual framework, the investigator examines the effect of $\bar{A}(K)$ on $Y(t)$ for all $t \in \mathcal{T}$ based on MSMs for the full data associated with $O$ whereas in the $t$-specific counterfactual framework, the investigator can examine the effect of $\bar{A}^t(t-s+1,t)$ on $Y^t$ for $t \in \mathcal{T}_s$ based on MSMs for the $t$-specific full data associated with $O^t$. Figure 5 illustrates how the $t$-specific counterfactual framework can be viewed as a collection of conventional counterfactual sub-frameworks with distinct definition of the outcome, treatment and baseline covariates. These conventional counterfactual sub-frameworks differ from the conventional counterfactual framework in the sense that the treatment history is of size $s \neq K + 1$ and the outcome of interest is no longer time-dependent.

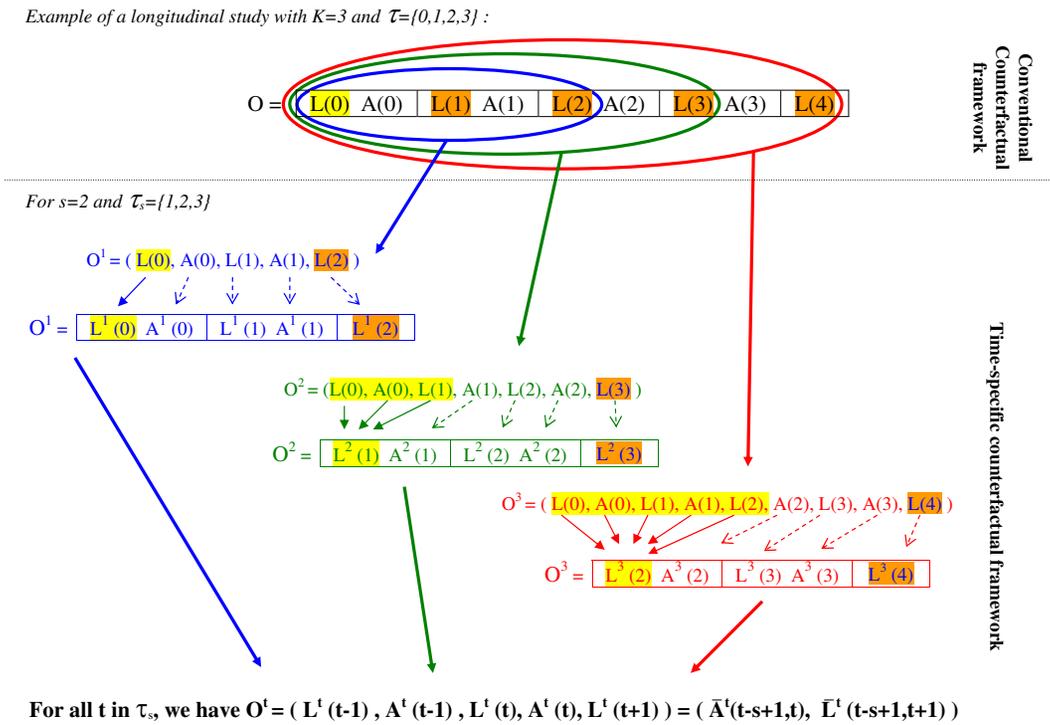

We have previously underscored that the $t$-specific and conventional counterfactual approaches lead to a different representation of the causal effect of $A$ on $Y$ with MSMs. We have also shown the equivalence between the parameters defined with MSMs in the $t$-specific counterfactual framework and HRMSM parameters defined in the conventional counterfactual framework and argued that these parameters may often provide a representation of causal effects which are more relevant for public health research. We can thus now easily address the issue of HRMSM parameter estimation. HRMSM estimators can indeed be derived from the application of the conventional MSM estimation methodologies to all the aforementioned conventional counterfactual sub-frameworks.

### A.5. HRMSM estimators

Under the assumptions presented earlier in section A.2 and from equality (6), the HRMSM parameters, $\beta_t(F_X, g \mid \cdot)$, can be identified and consistently estimated with the $t$-specific observed data and three estimators of the $t$-specific MSM parameters: the Inverse Probability of Treatment Weighted (16, 17, 19), the G-computation (12, 13, 3, 15, 29, 28, 8, 9) and Double Robust (18, 26, 7) estimators. The implementation procedures for these three estimators of HRMSM parameters correspond with the procedures developed for MSM-based causal inference except that they are applied not to the observed data, $O$, with treatment, $A$, and time-dependent outcomes of interest, $Y(t)$ for $t \in \mathcal{T}$, but to all $t$-specific observed data, $O^t$, of interest, i.e. for $t \in \mathcal{T}_s$, with treatment $A^t$ and outcome $Y^t$. The consistency and efficiency properties of these three estimators along with their implementation procedures have been thoroughly studied in the literature (see references above).

## Appendix B: Sufficient assumptions for consistent HRMSM estimation

We now formally establish that the seemingly larger set of $t$-specific assumptions (see section A.2) required to investigate causal effects in the $t$-specific counterfactual framework is implied by the set of assumptions (see section 2.2) required to investigate causal effects in the conventional counterfactual framework. The important practical consequence of both of these results is that successful investigation of causal effects with HRMSM parameters can be achieved based on the same model assumptions leading to successful investigation of causal effects with MSMs. Thus, the choice of HRMSM-based causal analysis over MSM-based causal analysis is only a matter of practical considerations (statistical power and computing issues) and above all subject-matter considerations (the relevance of the causal effect representation to public health research).

**Theorem B.1.** *We adopt the notations introduced previously for the conventional and $t$-specific counterfactual frameworks and in particular $t^- \equiv t - s + 1$. Based on these notations we have:*

  i. *the assumption of existence of counterfactuals defined in the conventional counterfactual framework implies the $t$-specific assumptions of existence of counterfactuals defined in the $t$-specific counterfactual framework:*

$$\forall \; \bar{a}(K) \in \mathcal{A}_V(K) \;\; \bar{L}_{\bar{a}(K)}(K+1) \;\; \implies \;\; \forall \; t \in \{s-1, \dots, K\} \; \forall \; \bar{a}^t(t^-, t) \in$$
$$\mathcal{A}_V(t^-, t) \;\; \bar{L}^t_{\bar{a}^t(t^-,t)}(t^-, t+1)$$

  ii. *the consistency assumption in the conventional counterfactual framework implies the $t$-specific consistency assumptions in the $t$-specific counterfactual framework:*

$$\forall \; t \in \{0, \dots, K+1\} \;\; L(t) = L_{\bar{A}(K)}(t) \;\; \implies \;\; \forall \; t \in \{s-1, \dots, K\} \;\; \forall \; j \in$$
$$\{t^-, \dots, t+1\} \quad L^t(j) = L^t_{\bar{A}^t(t^-,t)}(j)$$

  iii. *the temporal ordering assumption in the conventional counterfactual framework implies the $t$-specific temporal ordering assumptions in the $t$-specific counterfactual framework:*

$$\forall \; t \in \{0, \dots, K+1\} \quad L_{\bar{a}(K)}(t) = L_{\bar{a}(t-1)}(t) \;\; \implies \;\; \forall \; t \in \{s-1, \dots, K\}$$
$$\forall \; j \in \{t^-, \dots, t+1\} \quad L^t_{\bar{a}^t(t^-,t)}(j) = L^t_{\bar{a}^t(t^-,j-1)}(j)$$

*iv. the SRA in the conventional counterfactual framework implies the t-specific SRAs in the t-specific counterfactual framework:*

$$\forall\, t \in \{0, \ldots, K\} \quad A(t) \perp X \mid \bar{A}(t-1), \bar{L}(t) \implies \forall\, t \in \{s-1, \ldots, K\}$$

$$\forall\, j \in \{t^-, \ldots, t\} \quad A^t(j) \perp X^t \mid \bar{A}^t(t^-, j-1), \bar{L}^t(t^-, j)$$

**Proof.** For $t \in \{s-1, \ldots, K\}$ and $j \in \{t^-, \ldots, t+1\}$ we have:

*iv.*

$$
\begin{aligned}
X^t &= \left( V, \bar{L}^t_{\bar{a}^t(t^-,t)}(t^-, t+1) \right)_{\bar{a}^t(t^-,t) \in \mathcal{A}_V(t^-,t)} \\
&= \left( V, \bar{L}^t_{\bar{a}(t^-,t)}(t^-, t+1) \right)_{\bar{a}(t^-,t) \in \mathcal{A}_V(t^-,t)} \text{ from (1)} \\
&= \left( V, L^t_{\bar{a}(t^-,t)}(t^-), \bar{L}^t_{\bar{a}(t^-,t)}(t^-+1, t+1) \right)_{\bar{a}(t^-,t) \in \mathcal{A}_V(t^-,t)} \\
&= \left( V, \bar{A}(t-s), \bar{L}_{\bar{a}(t^-,t)}(t^-), \bar{L}_{\bar{a}(t^-,t)}(t^-+1, t+1) \right)_{\bar{a}(t^-,t) \in \mathcal{A}_V(t^-,t)} \\
&\qquad \text{from (2) and (3)} \\
&= \left( V, \bar{A}(t-s), \bar{L}_{\bar{a}(t^-,t)}(t+1) \right)_{\bar{a}(t^-,t) \in \mathcal{A}_V(t^-,t)} \\
&= \left( V, \bar{A}(t-s), \left( \bar{L}_{\bar{a}(t^-,t)}(t+1) \right)_{\bar{a}(t^-,t) \in \mathcal{A}_V(t^-,t)} \right) \\
X^t &= (V, \bar{A}(t-s), X^t_L) \text{ where } X^t_L = \left( V, \bar{L}_{\bar{a}(t^-,t)}(t+1) \right)_{\bar{a}(t^-,t) \in \mathcal{A}_V(t^-,t)} \qquad (7)
\end{aligned}
$$

In addition, we have:

$$
\begin{aligned}
X^t_L &= \left( V, \bar{L}_{\bar{a}(t^-,t)}(t+1) \right)_{\bar{a}(t^-,t) \in \mathcal{A}_V(t^-,t)} \\
&= \left( V, \bar{L}_{\bar{A}(t-s),\bar{a}(t^-,t),\bar{A}(t+1,K)}(t+1) \right)_{\bar{a}(t^-,t) \in \mathcal{A}_V(t^-,t)} \\
&\subset X = \left( V, \bar{L}_{\bar{a}(K)}(K+1) \right)_{\bar{a}(K) \in \mathcal{A}_V(K)} \\
X^t_L &\subset X \qquad (8)
\end{aligned}
$$

Based on these previous two results, we obtain:

$$g\left(A^t(j), X^t \mid \bar{A}^t(t^-, j-1), \bar{L}^t(t^-, j)\right)$$

$$= g\left(A(j), X^t \mid \bar{A}(t^-, j-1), L^t(t^-), \bar{L}^t(t^- + 1, j)\right) \text{ from (1)}$$

$$= g\left(A(j), X^t \mid \bar{A}(t^-, j-1), \bar{A}(t-s), \bar{L}(t^-), \bar{L}(t^- + 1, j)\right) \text{ from (2) and (3)}$$

$$= g\left(A(j), X^t \mid \bar{A}(j-1), \bar{L}(j)\right)$$

$$= g\left(A(j), \bar{A}(t-s), X_L^t \mid \bar{A}(j-1), \bar{L}(j)\right) \text{ from (7)}$$

$$= g\left(A(j), X_L^t \mid \bar{A}(j-1), \bar{L}(j)\right) \text{ since } \bar{A}(t-s) \subset \bar{A}(j-1)$$

$$= g\left(A(j) \mid \bar{A}(j-1), \bar{L}(j)\right) \text{ from the SRA and (8)}$$

$$= g\left(A^t(j) \mid \bar{A}^t(t^-, j-1), \bar{L}^t(t^-, j)\right) \text{ from (1), (2) and (3)}$$

This last equality is equivalent to $A^t(j) \perp X^t \mid \bar{A}^t(t^-, j-1), \bar{L}^t(t^-, j)$.

We also have:

- if $j \neq t^-$:

  i.
  $$\begin{aligned}
  L^t_{\bar{a}^t(t^-, t)}(j) &= L^t_{\bar{a}(t^-, t)}(j) \text{ from (1)} \\
  &= L_{\bar{a}(t^-, t)}(j) \text{ from (2)} \\
  &= L_{A(0), \ldots, A(t-s), \bar{a}(t^-, t), A(t+1), \ldots, A(K)}(j)
  \end{aligned}$$

  ii.
  $$\begin{aligned}
  L^t(j) &= L(j) \text{ from (2)} \\
  &= L_{\bar{A}(K)}(j) \text{ from the consistency assumption} \\
  &= L_{\bar{A}(t-s), \bar{A}(t^-, t), \bar{A}(t+1, K)}(j) \\
  &= L_{\bar{A}(t^-, t)}(j) \\
  &= L_{\bar{A}^t(t^-, t)}(j) \text{ from (1)} \\
  &= L^t_{\bar{A}^t(t^-, t)}(j) \text{ from (2)}
  \end{aligned}$$

iii.

$$
\begin{aligned}
L^t_{\bar{a}^t(t^-,t)}(j) &= L_{\bar{a}(t^-,t)}(j) \text{ from (1) and (2)} \\
&= L_{\bar{A}(t-s),\bar{a}(t^-,t),\bar{A}(t+1,K)}(j) \\
&= L_{\bar{A}(t-s),\bar{a}(t^-,j-1)}(j) \text{ from the temporal ordering} \\
&\quad \text{assumption} \\
&= L_{\bar{a}(t^-,j-1)}(j) \\
&= L^t_{\bar{a}^t(t^-,j-1)}(j) \text{ from (1) and (2)}
\end{aligned}
$$

- if $j = t^-$:

  i.

$$
\begin{aligned}
L^t_{\bar{a}^t(t^-,t)}(j) &= L^t_{\bar{a}(t^-,t)}(j) \text{ from (1)} \\
&= \left(\bar{A}(0,t-s), \bar{L}(0,t^-)\right)_{\bar{a}(t^-,t)} \text{ from (3)} \\
&= \left(\bar{A}(t-s), \bar{L}_{\bar{a}(t^-,t)}(t^-)\right) \\
&= \left(\bar{A}(t-s), \bar{L}_{A(0),\dots,A(t-s),\bar{a}(t^-,t),A(t+1),\dots,A(K)}(t^-)\right)
\end{aligned}
$$

  ii.

$$
\begin{aligned}
L^t(j) &= \left(\bar{A}(0,t-s), \bar{L}(0,t^-)\right) \text{ from (3)} \\
&= \left(\bar{A}(0,t-s), \bar{L}_{\bar{A}(K)}(0,t^-)\right) \text{ from the consistency} \\
&\quad \text{assumption} \\
&= \left(\bar{A}(0,t-s), \bar{L}_{\bar{A}(t-s),\bar{A}(t^-,t),\bar{A}(t+1,K)}(0,t^-)\right) \\
&= \left(\bar{A}(0,t-s), \bar{L}_{\bar{A}(t^-,t)}(0,t^-)\right) \\
&= \left(\bar{A}(0,t-s), \bar{L}(0,t^-)\right)_{\bar{A}(t^-,t)} \\
&= L^t_{\bar{A}(t^-,t)}(j) \text{ from (3)} \\
&= L^t_{\bar{A}^t(t^-,t)}(j) \text{ from (1)}
\end{aligned}
$$

iii.

$$
\begin{aligned}
L^t_{\bar{a}^t(t^-,t)}(j) &= \left( \bar{A}(0,t-s), \bar{L}(0,t^-) \right)_{\bar{a}(t^-,t)} \text{ from (1) and (3)} \\[2mm]
&= \left( \bar{A}(0,t-s), \bar{L}_{\bar{a}(t^-,t)}(0,t^-) \right) \\[2mm]
&= \left( \bar{A}(0,t-s), \bar{L}_{\bar{A}(t-s),\bar{a}(t^-,t),\bar{A}(t+1,K)}(0,t^-) \right) \\[2mm]
&= \left( \bar{A}(0,t-s), \bar{L}_{\bar{A}(t-s)}(0,t^-) \right) \text{ from the temporal} \\
&\qquad \text{ordering assumption} \\[2mm]
&= \left( \bar{A}(0,t-s), \bar{L}(0,t^-) \right)_{\bar{A}(t-s)} \\[2mm]
&= L^t_{\bar{A}(t-s)}(j) \text{ from (3)} \\[2mm]
&= L^t_{\bar{A}(t-s),\bar{a}^t(t^-,j-1)}(j) \text{ since } \bar{a}^t(t^-,t-s) \text{ is empty by} \\
&\qquad \text{definition} \\[2mm]
&= L^t_{\bar{a}^t(t^-,j-1)}(j) \square
\end{aligned}
$$

**Theorem B.2.** *We adopt the notations introduced previously for the conventional and t-specific counterfactual frameworks and in particular $t^- \equiv t - s + 1$. Based on these notations we have:*

*i.* $\forall\, t \in \{s-1,\dots,K\} \quad \forall\, j \in \{t^-,\dots,t+1\}$
$f(L^t(j) \mid \bar{L}^t(t^-,j-1), \bar{A}^t(t^-,j-1)) = f(L(j) \mid \bar{L}(j-1), \bar{A}(j-1))$

*ii.* $\forall\, t \in \{s-1,\dots,K\} \quad \forall\, j \in \{t^-,\dots,t\}$
$g(A^t(j) \mid \bar{A}^t(t^-,j-1), \bar{L}^t(t^-,j)) = g(A(j) \mid \bar{A}(j-1), \bar{L}(j))$

**Proof.**
For $t \in \{s-1,\dots,K\}$ we have

i. for $j \in \{t^-+1,\dots,t+1\}$:

$$
\begin{aligned}
&f\left( L^t(j) \mid \bar{L}^t(t^-,j-1), \bar{A}^t(t^-,j-1) \right) \\[2mm]
&= f\left( L^t(j) \mid L^t(t^-), \bar{L}^t(t^-+1,j-1), \bar{A}^t(t^-,j-1) \right) \text{ from (1)} \\[2mm]
&= f\left( L(j) \mid L^t(t^-), \bar{L}(t^-+1,j-1), \bar{A}^t(t^-,j-1) \right) \text{ from (2)} \\[2mm]
&= f\left( L(j) \mid \bar{A}(t-s), \bar{L}(t^-), \bar{L}(t^-+1,j-1), \bar{A}^t(t^-,j-1) \right) \text{ from (3)} \\[2mm]
&= f\left( L(j) \mid \bar{L}(j-1), \bar{A}(j-1) \right)
\end{aligned}
$$

for $j = t^-$:

$$f\left( L^t(j) \mid \bar{L}^t(t^-, j-1), \bar{A}^t(t^-, j-1) \right)$$

$$= f\left( L^t(j) \mid L^t(t^-), \bar{L}^t(t^-+1, j-1), \bar{A}(t^-, j-1) \right) \text{ from (1)}$$

$$= f\left( L^t(j) \mid L^t(t^-), \bar{L}(t^-+1, j-1), \bar{A}(t^-, j-1) \right) \text{ from (2)}$$

$$= f\left( L^t(j) \mid \bar{A}(t-s), \bar{L}(t^-), \bar{L}(t^-+1, j-1), \bar{A}(t^-, j-1) \right) \text{ from (3)}$$

$$= f\left( \bar{A}(j-1), \bar{L}(j) \mid \bar{L}(j-1), \bar{A}(j-1) \right) \text{ from (3)}$$

$$= f\left( L(j) \mid \bar{L}(j-1), \bar{A}(j-1) \right)$$

ii. for $j \in \{t^-, \ldots, t\}$:

$$g\left( A^t(j) \mid \bar{A}^t(t^-, j-1), \bar{L}^t(t^-, j) \right)$$

$$= g\left( A(j) \mid \bar{A}(t^-, j-1), L^t(t^-), \bar{L}^t(t^-+1, j) \right) \text{ from (1)}$$

$$= g\left( A(j) \mid \bar{A}(t^-, j-1), L^t(t^-), \bar{L}(t^-+1, j) \right) \text{ from (2)}$$

$$= g\left( A(j) \mid \bar{A}(t^-, j-1), \bar{A}(t-s), \bar{L}(t^-), \bar{L}(t^-+1, j) \right) \text{ from (3)}$$

$$= g\left( A(j) \mid \bar{A}(j-1), \bar{L}(j) \right)$$